\documentclass[17pt, fleqn, leqno]{article}
\usepackage{verbatim}
\usepackage[dvipsnames]{xcolor}
\usepackage[utf8]{inputenc}
\usepackage{extsizes}
\usepackage{graphicx}
\usepackage{amsmath}
\usepackage{amsfonts}
\usepackage{amssymb}
\usepackage{parskip}
\usepackage{setspace}
\usepackage{mathtools}
\usepackage{float}
\usepackage{dsfont}
\usepackage{tikz}
\usetikzlibrary{shapes,arrows,positioning,calc,intersections,patterns}
\usepackage [english]{babel}
\usepackage [autostyle, english = american]{csquotes}
\MakeOuterQuote{"}
\usepackage{pgfplots}
\pgfplotsset{compat=1.17}
\usepgfplotslibrary{fillbetween}
\usepackage[shortlabels]{enumitem}
\usepackage{hyperref}

\DeclareMathAlphabet{\mathbbold}{U}{bbold}{m}{n}
\makeatletter 
\def\@seccntformat#1{\@ifundefined{#1@cntformat} {\csname the#1\endcsname\quad} {\csname #1@cntformat\endcsname}}
\usepackage{titlesec}
\titleformat{\section}
{\large\normalfont}{\S\thesection.}{0.1cm}{}
\titleformat{\subsection}
{\large\normalfont}{\thesubsection}{0.1cm}{}
\title{\large\normalfont\color{Blue} On modular invariance of quantum affine $W$-algebras \vspace{-2ex}}
\author{Victor G. Kac and Minoru Wakimoto}
\date{\vspace{-4ex}}
\usepackage{cite}

\begin{document}

\maketitle

\underline{Abstract.} We find modular transformations of normalized characters for the following $W$-algebras:

\begin{enumerate}[(a)] 
    \item $W_k^{min}(\mathfrak{g}), \textnormal{where } \mathfrak{g}=D_n (n\geq4), \textnormal{or } E_6, E_7, E_8,$ and $k$ is a negative integer $\geq-2$, or $\geq -\frac{h^\vee}{6}-1$, respectively;
    \item quantum Hamiltonian reduction of the $\hat{\mathfrak{g}}$-module
 $L(k \Lambda_0)$, where $\mathfrak{g}$ is a simple Lie algebra, $f$ is its non-zero nilpotent element, and $k$ is a principal admissible level with the denominator $u>\theta(x)$, where $2x$ is the Dynkin characteristic of $f$ and $\theta$ is the highest root of $\mathfrak{g}$.
\end{enumerate}

We prove that these vertex algebras are modular invariant. A conformal vertex algebra $V$ is called modular invariant if  its character $tr_V q^{L_0-c/24}$ converges to a holomorphic modular function in the complex upper half-plane
on a congruence subgroup. We find explicit formulas for their characters.

Modular invariance of $V$ is important since, in particular, conjecturally it implies that
$V$ is simple, and that $V$ is rational, provided that it is lisse. 
\section{Introduction} 
It was proved in the paper \cite{KP84}, using the Weyl-Kac character formula, that the normalized characters of integrable irreducible highest weight modules of given level over an affine Lie algebra $\hat{\mathfrak{g}}$ over $\mathbb{C}$ are holomorphic modular functions, whose $\mathbb{C}$-span is $\textnormal{SL}_2(\mathbb{Z})$-invariant. In the remarkable paper \cite{Zhu96} this result was extended to an arbitrary rational vertex algebra $V$ for which all eigenvalues of $\textnormal{L}_0$ on $V$ are non-negative integers (some additional restrictions on $V$ in \cite{Zhu96} were removed in \cite{DLM98}.) In \cite{vEke13} Zhu's theorem was extended to $V$ with eigenvalues of $\textnormal{L}_0$ non-negative half-integers, which is the case for quantum affine $W$-algebras.

Next, in our paper \cite{KW88} we introduced the class of admissible highest weight $\hat{\mathfrak{g}}$-modules, which includes the integrable ones, and proved for them a generalization of the Weyl-Kac character formula. As for integrable modules in \cite{KP84}, the numerators of (normalized) characters of admissible modules are linear combinations of Jacobi theta forms, and since the (normalized) denominations are $\textnormal{SL}_2(\mathbb{Z})$-invariant, we derived in \cite{KW89}, Theorem 3.6, modular transformation formulas for these characters, extending the Kac-Peterson theorem to admissible modules.

We conjectured in \cite{KW89} that this modular invariance holds only for admissible $\hat{\mathfrak{g}}$-modules. This conjecture is still open.

Next, in \cite{KRW03} and \cite{KW04} we applied the quantum Hamiltonian reduction functor $H_f$ associated to a non-zero nilpotent element $f \in \mathfrak{g}$, to the universal affine vertex algebras $V^k(\mathfrak{g})$ of non-critical level k,
which produced the conformal vertex algebras denoted by $W^k(\mathfrak{g},f)$. We also considered the $W^k(\mathfrak{g},f)$-modules, obtained by applying $H_f$ to irreducible highest weight $\hat{\mathfrak{g}}$-modules $L(\Lambda)$; we write $H_f(\Lambda)$
for $H_f(L(\Lambda))$ to simplify notation. It is easy to see that the normalized Euler-Poincare character of $H_f(\Lambda)$ for an admissible $\hat{\mathfrak{g}}$-module $L(\Lambda)$ is still modular invariant. The simple quotient of
$W^k(\mathfrak{g},f)$
is denoted by $W_k(\mathfrak{g},f)$, and the latter vertex algebra for
$f=f_{min}=e_{-\theta}$ is denoted by $W_k^{min}(\mathfrak{g})$.

In \cite{Kaw18} Kawasetsu made a remarkable discovery that the $W$-algebra 
$W_k^{min}(\mathfrak{g})$  is a rational vertex algebra for the so called Deligne series $\mathfrak{g}=D_4, E_6,E_7, or E_8$, $k=-\frac{h^\vee}{6}$, 
%and $f$ an element from the non-zero nilpotent orbit of minimal dimension in $\mathfrak{g}$, 
hence modular invariance holds in spite of the fact that $-\frac{h^\vee}{6}\Lambda_0$ is not an admissible weight. (Recall that $W_{-\frac{h^\vee}{6}-1}^{min}(\mathfrak{g})$ is 1-dimensional \cite{AM18}.) Moreover, it was shown in \cite{AM18} that these $W$-algebras are lisse for negative integer levels $k \geq -\frac{h^\vee}{6}$, and the same holds if $\mathfrak{g}=D_n$ with $n>4$ for $k=-1, -2$. It was conjectured recently in \cite{ACK24} that these vertex algebras are rational, hence modular invariant.

In our paper \cite{KW18}, Theorem 4.1, we proved a character formula for $\hat{\mathfrak{g}}$-modules $L(\Lambda)$ of negative integer level $k$ and $\mathfrak{g}$ as above, under the assumptions that (i)$\Lambda$ is quasidominant, i.e.
$(\Lambda | \gamma) \in \mathbb{Z}_{\geq 0}$ for $\gamma \in\triangle_+$, (ii) there exists $\alpha \in \triangle_+$, for which $(\Lambda+\hat{\rho}|\alpha)=k+h^\vee$, (iii) if $\beta \in \hat{\triangle}_+$ is orthogonal to $\Lambda+\hat{\rho}$, then $\beta=\delta-\alpha$, and an extra hypothesis (iv), proved for some $\Lambda$ in \cite{BKK24}, including $\Lambda = k\Lambda_0$.

In the present note we prove, using these results, that the vertex algebra $W_k^{min}(\mathfrak{g})$ for the above $\mathfrak{g}$ and $k$ is modular invariant. If the (corrected) Dong-Ren conjecture in \cite{DR18}, that modular invariant lisse vertex algebra is rational, holds, then this result implies the ACK conjecture. We also prove that the vertex algebras
$\widetilde{W}_k(\mathfrak{g},f):=H_f(k\Lambda_0)$, where k is a principal admissible level with denominator $u$ and $f$ is a non-zero nilpotent element of
$\mathfrak{g}$ with Dynkin characteristic $2x$, such that $u>\theta(x)$, are modular invariant.

Note that the $W$-algebra $W_k(\mathfrak{g},f)$ is the simple quotient of
$\widetilde{W}_k(\mathfrak{g},f)$ (by Remark 5 below), and they coincide if the conjecture that modular invariance of a conformal vertex algebra implies its simplicity, holds.  

We keep the same notation of the theory of affine Lie algebras and their representations as in the book \cite{Kac90}.

\section{Main Theorem}

Let $\mathfrak{g}$ be a simple finite-dimensional Lie algebra over $\mathbb{C}$, let $\mathfrak{h} \subset \mathfrak{g}$ be a Cartan subalgebra, $\triangle \subset \mathfrak{h}^*$ the set of roots, and $\triangle_+$ a subset of positive roots. Let $\theta$ be the highest root and $(.|.)$ the invariant bilinear form on $\mathfrak{g}$, normalized by the condition $(\theta|\theta)=2$. Let $Q^\vee$ be the coroot lattice, and P the dual lattice. Let W be the Weyl group. 

Let $\hat{\mathfrak{g}}=\mathfrak{g}[t,t^{-1}]\oplus\mathbb{C}K\oplus\mathbb{C}d$ be the associated to $\mathfrak{g}$ affine Lie algebra and $\hat{\mathfrak{h}}=\mathfrak{h}\oplus\mathbb{C}K\oplus\mathbb{C}d$ its Cartan subalgebra. We coordinatize $\hat{\mathfrak{h}}$ by writing $h\in\hat{\mathfrak{h}}$ in the form $h=2\pi i( -\tau d+z+tK), \textnormal{where }\tau,t \in \mathbb{C} \textnormal{ and }z\in\mathfrak{h}$, and write a function $f$ on $\hat{\mathfrak{h}}$ in the form $f=f(\tau,z,t)$; we also write $f(\tau,z)=f(\tau,z,0)$. We identify $\mathfrak{h}$ with $\mathfrak{h}^*$ and $\hat{\mathfrak{h}}$ with $\hat{\mathfrak{h}}^*$ using $(.|.)$ and the usual extension of $(.|.)$ from $\mathfrak{g}$ to an invariant bilinear form on $\hat{\mathfrak{g}}$ \cite{Kac90}.

For $\mathfrak{g}$ of type $D_n(n\geq4), E_6,E_7, \textnormal{or }E_8 \textnormal{ let } b=2 \textnormal{ if } \mathfrak{g}=D_n, \textnormal{ and } b=\frac{h^\vee}{6}+1$ if $\mathfrak{g}=E_6,E_7,\textnormal{or } E_8$. (Note that $b$ is the length of the longest leg of the Dynkin diagram of $\hat{\mathfrak{g}}$.) Then for any positive integer $j\leq b$ there exist unique simple roots $\alpha_1,...,\alpha_{j-1}$ of $\mathfrak{g}$, such that $\alpha^{(j)}\coloneqq\theta-\sum_{i=1}^{j-1}\alpha_i$ is a root. Since $(\rho|\theta)=h^\vee-1$, $(\rho|\alpha^{(-k)})=h^\vee+k$ for negative integers $k\geq-b$.

The conditions (i) - (iii) of \cite{KW18}, Theorem 4.1, hold for $\Lambda=k\Lambda_0$ and $\alpha=\alpha^{(-k)}$, and since the extra hypothesis (iv) holds as well \cite{BKK24}, we have the following character formula for $\Lambda=k\Lambda_0$ and some other $\Lambda \in \hat{\mathfrak{h}}^* \ \textnormal{and}\ \alpha$, listed there:
\begin{equation}
    \hat{R} \ \textnormal{ch}\ L(\Lambda)=\frac{1}{2}\sum_{w\in W}\epsilon(w)\sum_{\gamma\in Q^\vee}(\alpha|\gamma)e^{wt_\gamma(\Lambda+\hat{\rho})}.
\end{equation}
This formula is a motivation to introduce the following notation for $\Lambda\in \hat{\mathfrak{h}}^*$ of positive integer level $n=\Lambda(K)$:
\begin{equation}
    F_\Lambda^{[\alpha]}=\frac{1}{2}q^{|\Lambda|^2/2n}\sum_{\gamma\in Q^\vee}(\alpha|\gamma)e^{t_\gamma(\Lambda)},
  \end{equation}
where $q=e^{2\pi i\tau}$,  
\begin{equation}
    A_\Lambda^{[\alpha]}=\sum_{w\in W}\epsilon(w)wF_\Lambda^{[\alpha]}.
\end{equation}

% \section{}
Let $\{e,x,f\}$ be an $\textnormal{sl}_2$-triple in $\mathfrak{g}$, so that $[e,f]=x,\ [x,e]=e,\ [x,f]=-f$, $x\in \mathfrak{h},\ e\in \mathfrak{n}_{+},\ f\in \mathfrak{n}_{-}$, and we have the $ad\ x$-eigenspace decomposition, compatible with the triangular decomposition:
\begin{equation}
    \mathfrak{g}=\bigoplus_{j\in\frac{1}{2}\mathbb{Z}}\mathfrak{g}_j.
\end{equation}
Since $e=\sum_{\beta\in B}e_{\beta}$ for some subset $B \subset\triangle_{+}$, and $[e,\mathfrak{h}^f]=0$, we see that $\beta|_{\mathfrak{h}^f}=0$ for all $\beta\in B$. Since $\mathfrak{g}_1\neq0$, there exists $\beta\in B$, for which $\beta(x)=1$. We fix $\beta\in \triangle_+$, such that
\begin{equation}
    \beta|_{\mathfrak{h}^f}=0\ \textnormal{and}\ \beta(x)\in\mathbb{Z}.
\end{equation}
The following is the key result of the paper.

\underline{Theorem 1}. Let $\Lambda\in \hat{\mathfrak{h}}^*$ be an integral weight for $\hat{\mathfrak{g}}$ of positive level $n=\Lambda(K)$, and let $\alpha\in\triangle_+$ be such that $(\Lambda|\alpha)=n$. Choose $\beta\in\triangle_+$, satisfying (5). Then for $z\in\mathfrak{h}^f$ one has
\begin{equation}
    \begin{split}
        &A_\Lambda^{[\alpha]}(\tau, -\tau x+z, \frac{1}{2}\tau|x|^2)\\ 
        =&\frac{1}{4}\beta(x)\sum_{w\in W}\epsilon(w)(w(\alpha)|\beta^\vee)\sum_{\gamma\in Q^\vee}e^{2\pi i(\overline{\Lambda}+n\gamma|w^{-1}(z))}q^{|\overline{\Lambda}+n(\gamma-w^{-1}(x))|^2/2n},
    \end{split}
\end{equation}
where $\overline{\Lambda}=\Lambda|_\mathfrak{h}$.

The proof is based on two lemmas.

\underline{Lemma 1}. For any $w\in W$ and $z\in\mathfrak{h}^f$ the following formula holds:
\begin{equation}
    \begin{split}
        & (wF_\Lambda^{[\alpha]}-r_\beta wF_\Lambda^{[\alpha]})(\tau,-\tau x+z, \frac{\tau}{2}|x|^2)\\
        = & \frac{\beta(x)}{2}(w\alpha|\beta^\vee)\sum_{\gamma\in Q^\vee}e^{2\pi i(\overline{\Lambda}+n\gamma|w^{-1}z)}q^\frac{|\overline{\Lambda}+n(\gamma-w^{-1}x)|^2}{2n}.
    \end{split}
\end{equation}
\underline{Proof}. Since
\begin{equation*}
    wF_\Lambda^{[\alpha]}=\frac{1}{2}\sum_{\gamma\in Q^\vee}(\gamma|\alpha)e^{w(\Lambda+n\gamma)}q^\frac{|\overline{\Lambda}+n\gamma|^2}{2n},
\end{equation*}
we have
\begin{equation}
    (wF_\Lambda^{[\alpha]})(\tau,-\tau x+z, \frac{\tau}{2}|x|^2)=\frac{1}{2}\sum_{\gamma\in Q^\vee}(\gamma|\alpha)e^{2\pi i(\overline{\Lambda}+n\gamma|w^{-1}z)}q^{|\overline{\Lambda}+n(\gamma-w^{-1}x)|^2/2n},
\end{equation}
\begin{equation}
    \begin{split}
        & (r_\beta wF_\Lambda^{[\alpha]})(\tau,-\tau x+z, \frac{\tau}{2}|x|^2)=\\
        & \frac{1}{2}\sum_{\gamma\in Q^\vee}(\gamma|\alpha)e^{2\pi i(\overline{\Lambda}+n\gamma|w^{-1}z)}q^{|\overline{\Lambda}+n(\gamma-w^{-1}x)+n(x|\beta)w^{-1}\beta^\vee|^2/2n}.
    \end{split}
\end{equation}
For (8) we used that $r_\beta(z)=z$ by (5).
Letting $\gamma'=\gamma+\beta(x)w^{-1}(\beta^\vee)$, and using again that $\beta(z)=0$, we rewrite (9) as (here we use that $\beta(x)\in\mathbb{Z}$)
\begin{equation*}
    \begin{split}
        & (r_\beta wF_\Lambda^{[\alpha]})(\tau,-\tau x+z, \frac{\tau}{2}|x|^2)\\
        = & \frac{1}{2}\sum_{\gamma'\in Q^\vee}(\gamma'-\beta(x)w^{-1}\beta^\vee|\alpha)e^{2\pi i(\overline{\Lambda}+n\gamma'|w^{-1}z)}q^\frac{|\overline{\Lambda}+n(\gamma'-w^{-1}x)|^2}{2n}.
    \end{split}
\end{equation*}
Subtracting this equality from (8), we obtain (7).

\underline{Lemma 2}. Let $W^\beta$ be a subset of $W$, such that 
    $W=W^\beta \coprod r_\beta W^\beta$.
Then 
\begin{equation}
    \begin{split}
        & A_\Lambda^{[\alpha]}(\tau,-\tau x+z, \frac{\tau}{2}|x|^2)\\
        = & \frac{1}{2}\beta(x)\sum_{w\in W^\beta}\epsilon(w)(w(\alpha)|\beta^\vee)\sum_{\gamma\in Q^\vee}e^{2\pi i(\overline{\Lambda}+n\gamma|w^{-1}z)}q^{|\overline{\Lambda}+n(\gamma-w^{-1}x)|^2/2n}.
    \end{split}
\end{equation}
\underline{Proof}. We have
\begin{equation*}
    A_\Lambda^{[\alpha]}(\tau,-\tau x+z, \frac{\tau}{2}|x|^2) = \sum_{w\in W^\beta}\epsilon(w)(wF_\Lambda^{[\alpha]}-r_\beta wF_\Lambda^{[\alpha]})(\tau,-\tau x+z, \frac{\tau}{2}|x|^2).
\end{equation*}
Using Lemma1, Lemma 2 follows.

\underline{Proof of Theorem 1}. Note that we can replace $W^\beta$ by $r_\beta W^\beta$ in Lemma 2, obtaining
\begin{equation*}
    \begin{split}
        & A_\Lambda^{[\alpha]}(\tau,-\tau x+z, \frac{\tau}{2}|x|^2)\\
        = &\frac{1}{2}\beta(x)\sum_{w\in r_\beta W^\beta}\epsilon(w)(w(\alpha)|\beta^\vee)\sum_{\gamma\in Q^\vee}e^{2\pi i(\overline{\Lambda}+n\gamma|w^{-1}
          (z))}q^{|\overline{\Lambda}+n(\gamma-w^{-1}x)|^2/2n}.
    \end{split}
\end{equation*}
Adding (10) to this equality and dividing by 2, we obtain (6).

Obviously (6) is zero if $f$ is a distinguished nilpotent element, different from the principal one. Indeed, in this case
$\mathfrak{h}^f=0$  and there exists $\beta\in\triangle $, such that $\beta(x)=0$.
A more detailed analysis in \cite{Wak24} shows that (6) is zero, unless $f$ lies in the adjoint orbit of $e_{-\theta}$ (=non-zero adjoint orbit of minimal dimension), when $\mathfrak{g}$ is simply laced. In this case there is a unique $\beta\in\triangle_+$, satisfying (5), namely, $\beta=\theta$.

\underline{Remark 1}. 
For $\gamma \in \Delta_+$, define the derivation $D_{\gamma}$ by 
$D_{\gamma}(e^{\mu}) = (\gamma|\mu) e^{\mu}$, and let 
$D = \prod_{\gamma \in \Delta_+}D_{\gamma}$.  Applying $D$ to both sides of 
(1) for $\Lambda =k \Lambda_0$, and letting after that 
$z=0$, $t=0$, we obtain for $\frak{g}$, $k$ and $\alpha$ as above,
$$
tr_{V_k(\mathfrak{g})}q^{L_0-\frac{c(\mathfrak{g},k)}{24}}
= 
\sum_{\gamma \in Q^{\vee}}(\alpha|\gamma) d((k+h^{\vee})\gamma)
q^{\frac{|\rho+(k+h^{\vee})\gamma|^2}{2(k+h^{\vee})}}/
2\eta(\tau)^{{\rm dim} \mathfrak{g}} ,
$$
where $c(\mathfrak{g},k)$ and $d(\gamma)$ are as in Example 2 in Section 4.

If $\Lambda$ is quasi-dominant of level $k > -h^{\vee}$, then
$$
ch L(\Lambda) = \sum_{w \in W}\epsilon(w) \sum_{\gamma \in Q^{\vee}}
c(\gamma) e^{wt_{\gamma}(\Lambda+\hat{\rho})},
$$
where $c(\gamma) \in \mathbf{Z}$. Hence at $z=0$, $t=0$, we get by the 
same method:
$$
tr_{L(\Lambda)}q^{L_0-\frac{c(\mathfrak{g},k)}{24}}
= 
\sum_{\gamma \in Q^{\vee}}c(\gamma) 
d(\bar{\Lambda}+(k+h^{\vee})\gamma)
q^{\frac{|\bar{\Lambda}+\rho+(k+h^{\vee})\gamma|^2}{2(k+h^{\vee})}}/
\eta(\tau)^{{\rm dim} \mathfrak{g}} .
$$

\section{Jacobi forms and modular invariance of $W^{min}_k(\mathfrak{g})$ for negative integer $k$}

Formula (6) can be rewritten in terms of Jacobi forms of degree $n$ and rank $l$, which are analytic functions in the domain Im$\tau>0$, $z\in\mathfrak{h}^f$, (see e.g. \cite{Kac90}, Chapter 13 for a theory of Jacobi forms):
\begin{equation}
    \Theta_\lambda(\tau,z)=\sum_{\gamma\in Q^\vee}e^{2\pi i(\overline{\lambda}+n\gamma|z)}q^{|\overline{\lambda}+n\gamma|^2/2n},
\end{equation}
where $\lambda\in\hat{P}^n$, the set of integral weights of level $n>0$ for $\hat{\mathfrak{g}}$,
as follows:
\begin{equation}
    \begin{split}
        & A_\lambda^{[\alpha]}(\tau,-\tau x+z, \frac{\tau}{2}|x|^2)\\
        = & \frac{1}{4}\beta(x)\sum_{w\in W}\epsilon(w)(w\alpha|\beta^\vee)q^{n|x|^2/2}\Theta_\lambda(\tau,w^{-1}(z-\tau x)).
    \end{split}
\end{equation}
This motivates the introduction of the following functions:
\begin{equation}
    f_{\lambda,w}(\tau,z)=q^{n|x|^2/2}\Theta_\lambda(\tau,w^{-1}(z-\tau x)),
\end{equation}
\begin{equation}
    f^-_{\lambda,w}(\tau,z)=f_{\lambda,w}(\tau,z+x)=e^{4\pi i\lambda(x)}f_{\lambda,w}(\tau,z-x),
\end{equation}
\begin{equation}
    f^*_{\lambda,w}(\tau,z)=q^{-n|x|^2/2}f_{\lambda,w}(\tau,z+x+\tau x)=e^{4\pi i\lambda(x)}q^{-n|x|^2/2}f_{\lambda,w}(\tau,z-x+\tau x).
\end{equation}
The last two equalities in (14) and (15) hold since
\begin{equation}
    2\alpha(x)\in\mathbb{Z} \ \textnormal{for}\  \alpha\in\triangle.
\end{equation}

Recall that Jacobi forms satisfy the following modular transformation properties \cite{Kac90}, Theorem 3.5:
\begin{equation}
    \Theta_\lambda(-\frac{1}{\tau},\frac{z}{\tau})=\varphi(\tau,z)\sum_{\mu\in\hat{P}^n\mod nQ^\vee+\mathbb{C}\delta}a(\lambda,\mu)\Theta_\mu(\tau,z),
\end{equation}
where we use the following notations:
\begin{equation}
    \varphi(\tau,z)=(-i\tau)^{l/2}e^{\frac{\pi in}{\tau}|z|^2}|P/nQ^\vee|^{-\frac{1}{2}},\ l=\dim\mathfrak{h},
\end{equation}
\begin{equation}
    a(\lambda,\mu)=e^{-\frac{2\pi i}{n}(\overline{\lambda}|\overline{\mu})};
\end{equation}
\begin{equation}
    \Theta_\lambda(\tau+1,z)=e^{\frac{\pi i}{n}|\overline{\lambda}|^2}\Theta_\lambda(\tau,z).
\end{equation}

From these we deduce the following transformations of the functions (13)-(15), where $z\in\mathfrak{h}^f$ and $\mu$ runs over the set $\hat{P}^n\mod(nQ^\vee+\mathbb{C}\delta)$:
\begin{equation}
    f_{\lambda,w}(-\frac{1}{\tau},\frac{z}{\tau})=\varphi(\tau,z)\sum_\mu a(\lambda,\mu)f^*_{\mu,w}(\tau,z),
\end{equation}
\begin{equation}
    f_{\lambda,w}(\tau+1,z)=e^{-4\pi i\lambda(x)}e^{\pi i(n|x|^2+\frac{1}{n}|\overline{\lambda}|^2)}f^-_{\lambda,w}(\tau,z),
\end{equation}
\begin{equation}
    f^-_{\lambda,w}(-\frac{1}{\tau},\frac{z}{\tau})=\varphi(\tau,z)e^{4\pi i\lambda(x)}e^{-2\pi in|x|^2}\sum_\mu a(\lambda,\mu)f^-_{\mu,w}(\tau,z),
\end{equation}
\begin{equation}
    f^-_{\lambda,w}(\tau+1,z)=e^{\pi i(n|x|^2+\frac{1}{n}|\overline{\lambda}|^2)}f_{\lambda,w}(\tau,z),
\end{equation}
\begin{equation}
    f^*_{\lambda,w}(-\frac{1}{\tau},\frac{z}{\tau})=\varphi(\tau,z)e^{4\pi i\lambda(x)}\sum_\mu a(\lambda,\mu)f_{\mu,w}(\tau,z),
\end{equation}
\begin{equation}
    f^*_{\lambda,w}(\tau+1,z)=e^{\frac{\pi i}{n}|\overline{\lambda}|^2}f^*_{\lambda,w}(\tau,z).
\end{equation}

In order to construct an $\textnormal{SL}_2(\mathbb{Z})$-invariant family of functions, containing the characters of the QHR of the above series of vertex algebras, recall the following "numerator" functions, given by (12), which can be written as
\begin{equation}
    B_\lambda^{[\alpha]}(\tau, z) = \frac{1}{4}\beta(x)\sum_{w\in W}\epsilon(w)(w\alpha|\beta^\vee)f_{\lambda,w}(\tau,z),
\end{equation}
and introduce the functions $B_\lambda^{[\alpha]-}$ and $B_\lambda^{[\alpha]*}$, obtained from functions (27) replacing $f_{\lambda,w}$ by $f^-_{\lambda,w}$ and $f^*_{\lambda,w}$ respectively.

Recall also the normalized $W$-algebra denominator, obtained from the usual one in Remark 3.1 of \cite{KRW03}, multiplying by the factor
$i^{|\triangle_+^0|}e^{2\pi i\rho_0}q^{\frac{1}{16}\dim\mathfrak{g}_0-\frac{l}{48}}$ \cite{KW17}:
\begin{equation}
    \overset{W}{R}_f(\tau,z)=\eta(t)^{\frac{3}{2}l-\frac{1}{2}\dim\mathfrak{g}^f}\prod_{\alpha\in\triangle_+^0}
    \vartheta_{11}(\tau,\alpha(z))(\prod_{\alpha\in\triangle^{1/2}}\vartheta_{01}(\tau,\alpha(z)))^{1/2},
    \end{equation}
where $ z\in\mathfrak{h}^f$, $\triangle_+^0$ are positive roots in $\mathfrak{g}_0$, $\triangle^j$ are roots in $\mathfrak{g}_j$, and $\vartheta_{ab}(\tau,z)$, $a,b=0\ \textnormal{or}\ 1$, are the well-known Jacobi forms of rank 1 and degree 2 (see, e.g. \cite{KW14}, Appendix). In particular,
\begin{equation*}
    \vartheta_{11}(\tau,z)=-\vartheta_{11}(\tau,-z)=-iq^{\frac{1}{12}}e^{-\pi iz}\eta(\tau)\prod_{n=1}^\infty(1-e^{-2\pi iz}q^n)(1-e^{2\pi iz}q^{n-1}).
\end{equation*}
The other two denominators $\overset{W}{R_f^-}(\tau,z)$ and $\overset{W}{R_f^*}(\tau,z)$ are obtained from (28) replacing $z$ by $z+x$ and $z+x+\tau x$, respectively, and $\vartheta_{01}$ by $\vartheta_{00}$ and $\vartheta_{10}$, respectively.

The modular transformations of these functions follow from those of the $\vartheta_{ab}(\tau,z)$ (see, e.g. [KW14], Appendix):
\begin{equation}
    \overset{W}{R}_f(-\frac{1}{\tau},\frac{z}{\tau})=(-i)^{|\triangle_+^0|}(-i\tau)^{l/2}e^{\frac{\pi i}{\tau}A(z)}\overset{W}{R_f^*}(\tau,z),
\end{equation}
where
\begin{equation}
    A(z)=\frac{1}{2}\sum_{\alpha\in\triangle_0\cup\triangle_{1/2}}\alpha(z)^2,
\end{equation}
and similar formulas for $\overset{W}{R_f^-}$ and $\overset{W}{R_f^*}$, where
$\overset{W}{R_f^*}$ in the RHS of (29) is replaced by $\overset{W}{R_f^-}$ and $\overset{W}{R_f}$ respectively. Also we have
\begin{equation}
    \overset{W}{R_f}(\tau+1,z)=e^{\frac{\pi i}{12}(\dim\mathfrak{g}_0-\frac{1}{2}\dim\mathfrak{g}_{1/2})}\overset{W}{R_f^-}(\tau,z),
\end{equation}
\begin{equation}
    \overset{W}{R_f^-}(\tau+1,z)=e^{\frac{\pi i}{12}(\dim\mathfrak{g}_0-\frac{1}{2}\dim\mathfrak{g}_{1/2})}\overset{W}{R_f}(\tau,z),
\end{equation}
\begin{equation}
    \overset{W}{R^*}(\tau+1,z)=e^{\frac{\pi i}{12}\dim\mathfrak{g}^f}\overset{W}{R^*}(\tau,z).
\end{equation}
Along with the functions
\begin{equation}
    B_\lambda^{[\alpha]}(\tau, z) = \frac{1}{4}\beta(x)\sum_{w\in W}\epsilon(w)(w\alpha|\beta^\vee)f_{\lambda,w}(\tau,z),
\end{equation}
introduce the functions $B_\lambda^{[\alpha]-}$ and $B_\lambda^{[\alpha]*}$, replacing $f_{\lambda,w}$ by $f_{\lambda,w}^-$ and $f_{\lambda,w}^*$, respectively. Consider the functions
\begin{equation}
    \Psi_\lambda^{[\alpha]}(\tau,z)=\frac{B_\lambda^{[\alpha]}(\tau,z)}{\overset{W}{R_f}(\tau,z)},\Psi_\lambda^{[\alpha]-}(\tau,z)=\frac{B_\lambda^{[\alpha]-}(\tau,z)}{\overset{W}{R_f^-}(\tau,z)},\Psi_\lambda^{[\alpha]*}(\tau,z)=\frac{B_\lambda^{[\alpha]*}(\tau,z)}{\overset{W}{R_f^*}(\tau,z)}.
\end{equation}

From the modular transformations (20)-(26) of the functions $f_{\lambda,w}$, $f_{\lambda,w}^-$, $f_{\lambda,w}^*$ and (29)-(33) of the functions $\overset{W}{R_f}$, $\overset{W}{R_f^-}$, $\overset{W}{R_f^*}$, we obtain modular transformations of the functions (35):

\underline{Theorem 2}. Let $\alpha\in\triangle$, and let $\beta\in\triangle_+$ satisfy (5). Then for any $\lambda\in\hat{P}^n$ the functions (35) satisfy the following modular transformation properties ($\textnormal{Im}\tau>0, z\in\mathfrak{h}^f$):

(a)
\begin{equation*}
    \Psi_{\lambda}^{[\alpha]}(-\frac{1}{\tau}, \frac{z}{\tau}) = \varphi_{1}(\tau,z) \sum_{\mu}a(\lambda, \mu) \Psi_{\mu}^{[\alpha]*}(\tau, z),
\end{equation*}
\begin{equation*}
    \Psi_{\lambda}^{[\alpha]-}(-\frac{1}{\tau}, \frac{z}{\tau}) = \varphi_{1}(\tau,z) e^{2\pi i(2\lambda(x)-n|x|^{2})} \sum_{\mu}a(\lambda, \mu) \Psi_{\mu}^{[\alpha]-}(\tau, z),
\end{equation*}
\begin{equation*}
    \Psi_{\lambda}^{[\alpha]*}(-\frac{1}{\tau}, \frac{z}{\tau}) = \varphi_{1}(\tau,z) e^{4\pi i \lambda(x)} \sum_{\mu}a(\lambda, \mu) \Psi_{\mu}^{[\alpha]}(\tau, z),
\end{equation*}
where $a(\lambda, \mu)$ is given by (19), and (cf (18))
\begin{equation*}
    \varphi_{1}(\tau, z) = i^{|\triangle_{+}^{0}|} |P/nQ^\vee|^{-1/2} e^{\frac{\pi i}{\tau}(n|z|^{2}-\frac{1}{2}\sum\limits_{\gamma\in \triangle^{0}\bigcup\triangle^{1/2}} \gamma(z)^{2})}.
\end{equation*}
(b)
\begin{equation*}
    \Psi_{\lambda}^{[\alpha]}(\tau+1, z) = e^{\frac{\pi i}{n}|\overline{\lambda}-nx|^{2}} e^{-2\pi i (\lambda|x)} a_1 \Psi^{[\alpha]-}(\tau, z),
\end{equation*}
\begin{equation*}
    \Psi_{\lambda}^{[\alpha]-}(\tau+1, z) = e^{\frac{\pi i}{n}|\overline{\lambda}-nx|^{2}} e^{2\pi i (\lambda|x)} a_1 \Psi^{[\alpha]}(\tau, z),
\end{equation*}
\begin{equation*}
    \Psi_{\lambda}^{[\alpha]*}(\tau+1, z) = e^{\frac{\pi i}{n}|\overline{\lambda}|^{2} - \frac{\pi i}{12} \dim\mathfrak{g}^{f}}\Psi_{\lambda}^{[\alpha]*}(\tau, z),
\end{equation*}
where
\begin{equation}
    a_{1} = e^{-\frac{\pi i}{12}(\dim\mathfrak{g}_0-\frac{1}{2}\dim\mathfrak{g}_{1/2})}.
\end{equation}
\underline{Remark 2}. For $w'\in W$, we have:
\begin{equation*}
    f_{w'\lambda,w}(\tau, z) = f_{\lambda,ww'}(\tau,z),
\end{equation*}
hence $B_{w'\lambda}^{[w'\alpha]}(\tau,z)=\epsilon(w')B_{\lambda}^{[\alpha]}(\tau,z)$, and 
\begin{equation*}
    \Psi_{w'\lambda}^{[w'\alpha]}(\tau,z)=\epsilon(w')\Psi_{\lambda}^{[\alpha]}(\tau,z).
\end{equation*}
\underline{Remark 3}. In \cite{KW14}, formula (9.10), we introduced the normalized denominators $\overset{N}{R}\vphantom{R}_{(\epsilon')}^{(\epsilon)}(\tau, z)$, where $\epsilon, \epsilon'=0\ \textnormal{or}\ \frac{1}{2}$, for $N=2$ and
$N=4$ superconformal algebras. The functions $\overset{W}{R_f}$,
$\overset{W}{R^{-}_f}$, and $\overset{W}{R^{*}_f}$  
are analogues of $\overset{N}{R}\vphantom{R}_{(1/2)}^{(0)}$,
$\overset{N}{R}\vphantom{R}_{(1/2)}^{(1/2)}$, and
$\overset{N}{R}\vphantom{R}_{(0)}^{(1/2)}$, and the modular transformations, given by Proposition 9.1 there, are similar to (29)-(33).

%\section{}
Recall that the normalized character  $\textnormal{ch}_{\Lambda}$  of a $\hat{\mathfrak{g}}$-module $L(\Lambda)$ of level $k\neq-h^\vee$ and the normalized
Euler-Poincare character $\widetilde{\textnormal {ch}}_{H_f(\Lambda)}$ of the
$W^k(\mathfrak{g},f)$-module $H_f(\Lambda)$,
% obtained from $L(\Lambda)$ by the QHR, associated to $f$,
are related by the following formula (see \cite{KRW03}, and \cite{KW17},
formula (9)):
\begin{equation}
  (\overset{W}{R_f}\ \widetilde{\textnormal{ch}}_{H_f(\Lambda)})(\tau,z)=(\hat{R}\
  \textnormal{ch}_{\Lambda})
  (\tau,-\tau x+z, \frac{\tau}{2}|x|^2), z\in\mathfrak{h}^f.
\end{equation}
Recall also that if $\mathfrak{g}$ is of type $D_n (n\geq4),\ E_6,\ E_7,\ \textnormal{or}\ E_8$, $b$ is the positive integer introduced in \S 2, $k$ is a negative integer $\geq-b$, $\alpha=\alpha^{(-k)}$, then we have the character formula (1) for $\Lambda=k\Lambda_0$ (where $\hat{R}$ is normalized by the factor $e^{\hat{\rho}}$). Using \cite{Ara05}, $H_{f_{min}}(k \Lambda_0)=W_k^{min}(\mathfrak{g})$.
% the QHR, associated to $f_{min}:=e_{-\theta}$, of $\textnormal{L}
Consequently, by Theorem 1 and (37) we obtain
\begin{equation}
  \textnormal{ch}_{W_k^{min}(\mathfrak{g})}(\tau, z)=
  \textnormal{ch}_{H_{f_{min}}(k\Lambda_0)}(\tau,z)=\Psi^{[\alpha]}_{k\Lambda_0+\hat{\rho}}(\tau,z), \,z\in\mathfrak{h}^f.
\end{equation}
By \cite{AM18} the vertex algebra $W_k^{min}(\mathfrak{g})$ is lisse, hence by \cite{Zhu96}, the function $\Psi^{[\alpha]}_{k\Lambda_0}(\tau,z)$ is analytic in a neighborhood of $z=0$ in $\mathfrak{h}^f$ for each $\tau$ with  Im $\tau>0$.
(An elementary proof of this fact follows from the observation that if
$\gamma\in\triangle^0$, and $\gamma$ vanishes on $\mathfrak{h}^f$,
then $f_{\lambda, w}$ is unchanged if $w$ is replaced by $r_{\gamma}w$.)
By Theorem 2, this function lies in the finite-dimensional $\textnormal{SL}_2(\mathbb{Z})$-invariant space, which is the $\mathbb{C}$-span of the functions
\begin{equation*}
    \begin{split}
      \Psi_{\Lambda+\hat{\rho}}^{[\alpha]}(\tau,z),\, \Psi_{\Lambda+\hat{\rho}}^{[\alpha]-}(\tau,z),\,\Psi_{\Lambda+\hat{\rho}}^{[\alpha]*}(\tau,z),
      \end{split}
\end{equation*}
where
$\alpha\in\triangle_+,
\Lambda+\hat{\rho}\in\hat{P}^{k+h^\vee}\mod(k+h^\vee)Q^\vee+\mathbb{C}\delta$.
Thus we obtain the following theorem, using also Proposition 13.6 from \cite{Kac90}
for $\alpha=0, \beta=nw^{-1}(x)$ there:

\underline{Theorem 3}. For $\mathfrak{g}$ and $k$ as above, the normalized character of the vertex algebra $W_k^{min}(\mathfrak{g})$, is given by
\begin{equation}
    tr_{W_k^{min}(\mathfrak{g})}q^{L_0-c(k)/24}=\lim_{z\rightarrow 0} \Psi_{k\Lambda_0+\hat{\rho}}^{[\alpha]}(\tau,z),
\end{equation}
and it is a holomorphic modular function with respect to a congruence subgroup, on the upper half-plane Im$\tau>0$ ($c(k)$ is given by \cite{KW04}, formula (5.7)).

\underline{Remark 4}. Let $\mathfrak{g}$ be one of the Lie algebras $D_4,E_6,E_7,E_8$. Then $\dim W^b_{min}(\mathfrak{g})=1$ \cite{AM18}, hence we obtain from (38), (35) the following denominator identity for the normalized denominator of $W_k^{min}(\mathfrak{g})$, given by (28) for $f=e_{-\theta}$:
\begin{equation*}
    \begin{split}
      &\overset{W}{R_f}(\tau,z)=\frac{1}{4}\sum_{w\in W}\epsilon(w)
      (w(\alpha)|\theta)\\ 
        & \times\sum_{\gamma\in Q^\vee}e^{2\pi i(\rho+(h^\vee-b)\gamma|w^{-1}(z))}q^{|\rho+(h^\vee-b)(\gamma-w(\theta)/2)|^2/2(h^\vee-b)}.
    \end{split}
\end{equation*}

\section{Modular invariance at principal admissible levels}
Recall that the principal admissible weights for $\hat{\mathfrak{g}}$ \cite{KW89} is the set, denoted $\hat{Pr}^k$, consisting of weights $\lambda$ of
principal admissible level $k$, which is a rational number with denominator $u\in\mathbb{Z}_{\geq1}$, such that
\begin {equation}
k+h^\vee = \frac{p}{u},\, 
\gcd(p,u)=\gcd(u,r^\vee)=1,\, p\geq h^\vee,
\end{equation}
where $r^\vee$ is the ``lacety'' of $\mathfrak{g}$, and which are of the form
\begin {equation}
\lambda=(t_\beta y). (\Lambda^0-(u-1)(k+h^\vee)\Lambda_0),
\end{equation}
where $\beta\in Q^*, y\in W$ are such that $(t_\beta y)\{u\delta-\theta, \alpha_1,...,\alpha_l\}\subset\hat{\triangle}_+$, and $\Lambda^0$ is a dominant integral weight of level $u(k+h^\vee)-h^\vee$.

\underline{Remark 5}. By the result of Arakawa in \cite{Ara15}, for any $\mathfrak{g}$-locally finite $\hat{\mathfrak{g}}$-module $L(\lambda)$, where $\lambda$ is an admissible weight, and any non-zero nilpotent element $f$ of $\mathfrak{g}$, $H_f(\lambda)=H_f^0(\lambda)$. It follows that, in particular, the Euler-Poincare character of $H_f(\lambda)$ coincides with its character, and that if $H_f(k\Lambda_0)$ is not zero, then the simple $W$-algebra $W_k(\mathfrak{g},f)$ is a quotient of $H_f(k\Lambda_0)=\widetilde{W}_k(\mathfrak{g},f)$.
If the conjecture that any modular invariant conformal vertex algebra is simple, holds, then we may conclude that $H_f(k\Lambda_0)$ is either zero, or is isomorphic to $W_k(\mathfrak{g},f)$.

\underline{Remark 6}. Since a
$\hat{\mathfrak{g}}$-module $L(\lambda)$ is
$\mathfrak{g}$-locally finite iff
$\lambda$ is quasidominant,
the principal admissible $\hat{\mathfrak{g}}$-module $L(\lambda)$ is
$\mathfrak{g}$-locally finite 
if $\beta=0$ and $y=1$ in (41) (and these are all $\mathfrak{g}$-locally finite principal admissible $L(\lambda)$). In particular, the principal admissible $\hat{\mathfrak{g}}$-module $L(k\Lambda_0)$ is $\mathfrak{g}$-locally finite. It follows from \cite{KW17}, formula (5), and Remark 5 that for these $\hat{\mathfrak{g}}$-modules the normalized character $ch_\lambda(\tau,z,t)$ is meronorphic in $\tau,z,t$, with $\tau,t \in\mathbb{C},\textnormal{Im}\tau>0, z\in\mathfrak{h}$, and analytic in a neighborhood of $z=0$ in $\mathfrak{h}$.

The modular transformations of the normalized characters $\textnormal{ch}_{\lambda}$ of $L(\lambda)$ with $\lambda\in\hat{Pr}^k$ are as follows \cite{KW89}, Theorem 3.6:
\begin{equation}
    \begin{split}
        & (\hat{R}\ \textnormal{ch}_{\lambda})(-\frac{1}{\tau},\frac{z}{\tau},t-\frac{|z|^2}{2\tau})\\
        & =(-i\tau)^{l/2}(-i)^{|\triangle_+|}\sum_{\mu\in\hat{Pr}^k\mod\mathbb{C}\delta}a(\lambda,\mu)(\hat{R}\ \textnormal{ch}_{\mu})(\tau,z,t),
    \end{split}
\end{equation}
where the numbers $a(\lambda,\mu)$ are given in \cite{KW89}, Theorem 3.6;
\begin{equation}
    (\hat{R}\ \textnormal{ch}_{\lambda})(\tau+1,z,t)=e^{\frac{\pi i}{k+h^\vee}|\overline{\lambda}+\rho|^2}(\hat{R}\ \textnormal{ch}_{\lambda})(\tau,z,t).
\end{equation}

As in (37) and (13), (14), the normalized Euler-Poincare character, the signed one, and the signed twisted one $\widetilde{\textnormal{ch}}_{H_f(\lambda)}$, $\widetilde{\textnormal{ch}}_{H_f(\lambda)}^-$ and $\widetilde{\textnormal{ch}}_{H_f(\lambda)}^*$, respectively,
% of the QHR of the $\hat{\mathfrak{g}}$-module $L(\lambda)$, $\lambda\in\hat{Pr}^k$,
are given by the following formulas $(z\in\mathfrak{h}^f)$:
\begin{equation}
    (\overset{W}{R_f}\ \widetilde{\textnormal{ch}}_{H_f(\lambda)})(\tau,z)=(\hat{R}\ \textnormal{ch}_{\lambda})(\tau,-\tau x+z, \frac{\tau}{2}|x|^2), 
\end{equation}
\begin{equation}
    (\overset{W}{R_f^-}\ \widetilde{\textnormal{ch}}_{H_f(\lambda)}^-)(\tau,z)=(\hat{R}\ \textnormal{ch}_{\lambda})(\tau,-\tau x+x+z, \frac{\tau}{2}|x|^2), 
\end{equation}
\begin{equation}
    (\overset{W}{R_f^*}\ \widetilde{\textnormal{ch}}_{H_f(\lambda)}^*)(\tau,z)=(\hat{R}\ {\textnormal{ch}}_{\lambda})(\tau, x+z, 0).
\end{equation}
As in the case of Theorem 2, it is straightforward to derive from (29)-(33) and (40), (41), the modular transformation properties of the functions
$\widetilde{\textnormal{ch}}_{H_f(\lambda)}$, $\widetilde{\textnormal{ch}}_{H_f(\lambda)}^-$, and $\widetilde{\textnormal{ch}}_{H_f(\lambda)}^*$, $\lambda\in\hat{Pr}^k$:

\underline{Theorem 4}. For $\lambda\in\hat{Pr}^k$ and $z\in\mathfrak{h}^f$ we have:

(a)
\begin{equation*}
    \widetilde{\textnormal{ch}}_{H_f(\lambda)}(-\frac{1}{\tau},\frac{z}{\tau})=\phi_2(\tau,z)\sum_{\mu\in{Pr}^k\mod\mathbb{C}\delta}a(\lambda,\mu)\widetilde{\textnormal{ch}}_{H_f(\mu)}^*(\tau,z),
\end{equation*}
\begin{equation*}
    \begin{split}
      \widetilde{\textnormal{ch}}_{H_f(\lambda)}^-(-\frac{1}{\tau},\frac{z}{\tau})=
      & \phi_2(\tau,z)e^{4\pi i(\lambda+\rho)(x)-2\pi i(k+h^\vee)|x|^2}\\
    & \times\sum_{\mu\in{Pr}^k\mod\mathbb{C}\delta}a(\lambda,\mu)\widetilde{\textnormal{ch}}_{H_f(\mu)}^-(\tau,z),
    \end{split}
\end{equation*}
\begin{equation*}
    \begin{split}
        \widetilde{\textnormal{ch}}_{H_f(\lambda)}^*(-\frac{1}{\tau},\frac{z}{\tau})= &\phi_2(\tau,z)e^{4\pi i(\lambda+\rho)(x)}\\
        &\times\sum_{\mu\in{Pr}^k\mod\mathbb{C}\delta}a(\lambda,\mu)\widetilde{\textnormal{ch}}_{H_f(\mu)}(\tau,z),
    \end{split}
\end{equation*}
where $a(\lambda,\mu)$ are as in \cite{KW89}, Theorem 3.6, and
\begin{equation*}
    \phi_2(\tau,z)=(-i)^{\dim \mathfrak{g}_{>0}}e^{\frac{\pi i}{\tau}((k+h^\vee)|z|^{2}-\frac{1}{2}\sum\limits_{\gamma\in \triangle^{0}\bigcup\triangle^{1/2}} \gamma(z)^{2})}.
\end{equation*}
(b)
\begin{equation*}
  \widetilde{ \textnormal{ch}}_{H_f(\lambda)}(\tau+1,z)=a_1e^{\frac{\pi i}{k+h^\vee}|\overline{\lambda}+\rho-(k+h^\vee)x|^2-2\pi i(\lambda+\rho)(x)}
  \widetilde{\textnormal{ch}}_{H_f(\lambda)}^-(\tau,z),
\end{equation*}
\begin{equation*}
  \widetilde{\textnormal{ch}}_{H_f(\lambda)}^-(\tau+1,z)=a_1e^{\frac{\pi i}{k+h^\vee}|\overline{\lambda}+\rho-(k+h^\vee)x|^2+2\pi i(\lambda+\rho)(x)}
  \widetilde{\textnormal{ch}}_{H_f(\lambda)}(\tau,z),
\end{equation*}
\begin{equation*}
  \widetilde{\textnormal{ch}}_{H_f(\lambda)}^*(\tau+1,z)= e^{\frac{\pi i}{k+h^\vee}|\overline{\lambda}+\rho|^2}e^{-\frac{\pi i}{12}\dim\mathfrak{g}^f}
  \widetilde{\textnormal{ch}}_{H_f(\lambda)}^*(\tau,z),
\end{equation*}
where $a_1$ is given by (36).

\underline{Remark 7}. For the $W$-algebra $W_k^{min}(\mathfrak{g})$, where $k$ is an admissible level, the QHR $H(\lambda)=H_{f_{min}}(\lambda)$ of an admissible $\hat{\mathfrak{g}}$-module $L(\lambda)$ is either 0, which happens if $(\lambda+\rho|\alpha_0)$ is a positive integer, or is irreducible, and in this case the $W_k^{min}(\mathfrak{g})$-modules $H(\lambda)$ and $H(\lambda')$ are isomorphic iff either $\lambda'=\lambda \mod\mathbb{C}\delta$ or $\lambda'=\lambda+(\lambda+\rho|\alpha_0)\theta\mod\mathbb{C}\delta$; also,
$H(\lambda)=H^0(\lambda)$. This follows from \cite{KW04} and \cite{Ara05}.
Using this remark one should be able to write the modular transformation formulas given by Theorem 4, in a basis of characters $\textnormal{ch}_{H(\lambda)}$, signed characters $\textnormal{ch}_{H(\lambda)}^-$, and signed twisted characters $\textnormal{ch}_{H(\lambda)}^*$.

%Using Remarks 5 and 6,
% that the simple affine vertex algebras $V_k(\mathfrak{g})$ at admissible level $k$ are quarilisse \cite{AK18},
%we obtain the following theorem (cf. Theorem 3).

%\underline{Theorem 5}.(a) The normalized character of a simple affine vertex algebra $V_k(\mathfrak{g})$ at admissible level $k$ is a holomorphic modular function on a congruence subgroup in the upper half-plane, given by the formula
%\begin{equation}
%    tr_{V_k(\mathfrak{g})}q^{L_0-c/24}=\lim_{z\rightarrow 0}\textnormal{ch}_{k%\Lambda_0}(\tau,z), \,z\in\mathfrak{h},
%\end{equation}
%where $c=\frac{k\dim\mathfrak{g}}{k+h^\vee}$ is the central charge of $V_k(\ma%thfrak{g})$.

%(b) The same claim as (a) holds if we replace in (45), $V_k(\mathfrak{g})$ by its QHR, associated to $f$, and $\lim_{z\rightarrow 0}\textnormal{ch}_{k\Lambda_0}(\tau,z)$ by $\lim_{z\rightarrow 0}\textnormal{ch}_{k\Lambda_0}(\tau,-\tau x+z, \frac{\tau}{2}|x|^2), z\in\mathfrak{h}^f$, and the central charge is given by \cite{KRW03}, Remark 2.2.

\underline{Example 1}. Let $\mathfrak{g}$ be a simple Lie algebra, and let $k_u=h^\vee\frac{1-u}{u}$ be a boundary admissible level, where $u\in\mathbb{Z}_{\geq1}$ is coprime to $h^\vee$ and to $r^\vee$. Then, by \cite{KW17}, Proposition 2b, we have
\begin{equation}
    tr_{V_{k_u}(\mathfrak{g})}q^{L_0-\frac{c(\mathfrak{g},k_u)}{24}}=\textnormal{ch}_{k_u\Lambda_0}(\tau,0)=\left(\frac{\eta(u\tau)}{\eta(\tau)}\right)^{\dim\mathfrak{g}},
\end{equation}
where $c(\mathfrak{g},k_u)=(1-u)\dim\mathfrak{g}$.

By \cite{KW17}, formula (11), we have for any non-zero nilpotent element $f$ of $\mathfrak{g}$ with Dynkin characteristic $2x$, and the boundary level $k_u$:
\begin{equation}
    \begin{split}
        tr_{H_f(k \Lambda_0)}q^{L_0-c(\mathfrak{g},f,k_u)/24}=(-i)^{|\triangle_+|}\frac{\eta(u\tau)^{\frac{3}{2}\dim\mathfrak{g}_0-\frac{1}{2}\dim\mathfrak{g}}}{\eta(\tau)^{\dim\mathfrak{g}_0-\dim\mathfrak{g}_{1/2}}\eta(\frac{\tau}{2})^{\dim\mathfrak{g}_{1/2}}}\\
        \times q^{\frac{h^\vee|x|^2}{2u}}\prod_{j\in\frac{1}{2}\mathbb{Z}_{>0}}\vartheta_{11}(u\tau,-j\tau)^{\dim\mathfrak{g}_j},
    \end{split}
\end{equation}
where $c(\mathfrak{g},f,k_u)=\dim\mathfrak{g}_0-\frac{1}{2}\dim\mathfrak{g}_{1/2}-\frac{12u}{h^\vee}|\rho-\frac{h^\vee}{u}x|^2$.

Note that (48) is a modular function since $q^{\frac{j^2}{2u}}\vartheta_{11}(u\tau,-j\tau)$ and $\eta(\tau)$  are modular forms of weight $\frac{1}{2}$ on a congruence subgroup. 
Note also that the RHS of (48) can be rewritten as follows:
\begin{equation*}
    \begin{split}
        (-i)^{|\triangle_+^0|}q^{\frac{h^\vee}{2u}|x|^2-\rho(x)+\frac{u}{24}(\dim\mathfrak{g}-\dim\mathfrak{g}_0)}\frac{\eta(u\tau)^{\dim\mathfrak{g}_0}}{\eta(\tau)^{\dim\mathfrak{g}_0-\dim\mathfrak{g}_{1/2}}\eta(\frac{\tau}{2})^{\dim\mathfrak{g}_{1/2}}}\\
        \times\prod_{j\in\frac{1}{2}\mathbb{Z}_{>0}}(\prod_{n=1}^\infty(1-q^{un-j})(1-q^{un-(u-j)}))^{\dim\mathfrak{g}_j},
    \end{split}
\end{equation*}

It follows that (48) is 0 iff $\dim\mathfrak{g}_u\neq0$. For the minimal
nilpotent $f_{min}$ this happens only if $u=1$, but for the principal nilpotent $f_{pr}$ this happens iff $u<h$, the Coxeter number of $\mathfrak{g}$.
By Remark 5, it follows that $\tilde{W}_{k_u}(\mathfrak{g},f)$ is a modular invariant vertex algebra, provided that $u>\theta(x)$.

The following proposition gives a beautiful product formula for the QHR, associated to $f_{pr}$ of $\mathfrak{g}$-locally finite principal admissible $\hat{\mathfrak{g}}$-modules (including $\textnormal{L}(k\Lambda_0)$), of level $k=-h^\vee+\frac{p}{h}$.

\underline{Proposition 1}. Let $f_{pr}$ be a principal nilpotent element of $\mathfrak{g}$ and let $k=-h^\vee+\frac{p}{h}$, where $h$ is the Coxeter number of $\mathfrak{g}, p\in\mathbb{Z}_{\geq h^\vee}, \gcd(p,h)=1$, and let $\lambda=\lambda^0-(h-1)(k+h^\vee)\Lambda_0$, where $\lambda^0\in\hat{P}_+^{p-h^\vee}$. Then the character of the QHR, associated to $f_{pr}$, of the $\hat{\mathfrak{g}}$-module $\textnormal{L}(\lambda)$ is given by the following formula
\begin{equation}
    \textnormal{ch}_{H_{f_{pr}}(\lambda)}=q^{|\overline{\lambda^0}+\rho-(k+h^\vee)\rho^\vee|^2/2(k+h^\vee)}\prod_{\alpha\in\hat{\triangle}_+^\vee}(1-q^{(\lambda^0+\hat{\rho}|\alpha)})^{mult\ \alpha}/\eta(\tau)^l.
\end{equation}
\underline{Proof}. Since $x=\rho^\vee$ and $\mathfrak{h}^f=0$, we have, by formula (9) from \cite{KW17}, and Remark 5:
\begin{equation*}
    \textnormal{ch}_{H_{f_{pr}}(\lambda)}=\frac{(\hat{R}\textnormal{ch}_{L(\lambda)})(\tau,-\tau\rho^\vee, \tau|\rho^\vee|^2/2)}{\eta(\tau)^l}.
\end{equation*}
By (5) from \cite{KW17}, we obtain from this:   
\begin{equation*}
    \textnormal{ch}_{H_{f_{pr}}(\lambda)}=q^{(k+h^\vee)|\rho^\vee|^2/2}A_{\lambda^0+\hat{\rho}}(h\tau,-\tau\rho^\vee,0)/\eta(\tau)^l,
\end{equation*}
where $A_{\lambda^0+\hat{\rho}}$ is the numerator of $\textnormal{ch}\ \textnormal{L}(\lambda^0)$. But $A_{\lambda^0+\hat{\rho}}(h\tau,-\tau\rho^\vee,0)=A_{\lambda^0+\hat{\rho}}(-2\pi i\tau\hat{\rho}^\vee)$, which, by the Weyl-Macdonald-Kac denominator identity, is equal to
\begin{equation*}
    q^{-(\lambda^0+\hat{\rho}|\hat{\rho}^\vee)}\prod_{\alpha\in\hat{\triangle}_+^\vee}(1-q^{(\lambda^0+\rho|\alpha)})^{mult\ \alpha}.
\end{equation*}
This completes the proof of (50).

\underline{Example 2}. let $\mathfrak{g}$  be a simple Lie algebra and let $k=-h^\vee+\frac{p}{u}$ be a principal admissible level (40).
%where $u\in\mathbb{Z}_{\geq 1}, p\in\mathbb{Z}_{\geq h^\vee}$ are coprime integers.
Then, by formula (3.3) from \cite{KW89} we have, using (46):
\begin{equation}
    tr_{V_k(\mathfrak{g})}q^{L_0-\frac{c(\mathfrak{g},k)}{24}}=\textnormal{ch}_{(p-h^\vee)\Lambda_0}(u\tau)\left(\frac{\eta(u\tau)}{\eta(\tau)}\right)^{\dim\mathfrak{g}},
\end{equation}
where $c(\mathfrak{g},k)=\frac{k\dim\mathfrak{g}}{k+h^\vee}$.
By \cite{Kac90}, Exercise 12.24, we have
\begin{equation}
     \textnormal{ch}_{(p-h^\vee)\Lambda_0}(\tau)=\frac{\sum_{\gamma\in pQ^\vee}d(\gamma)q^{|\rho+\gamma|^2/2p}}{\eta(\tau)^{\dim\mathfrak{g}}},
\end{equation}
where $d(\gamma)=\prod_{\alpha\in\triangle_+}\frac{(\gamma+\rho|\alpha)}{(\rho|\alpha)}$.

\underline{Theorem 5}. Let $\lambda\in\hat{Pr}^k$ be a principal admissible weight of level $k=\frac{p}{u}-h^\vee$, where $p\geq h^\vee$, of the form $\lambda=\Lambda^0-(u-1)(k+h^\vee)\Lambda_0$, where $\Lambda^0\in\hat{P}_+^{p-h^\vee}$ is a dominant integral weight (cf. (39) and (40)). Then the normalized character of the QHR, associated to a non-zero nilpotent element $f$, of the $\hat{\mathfrak{g}}$-module $\textnormal{L}(\lambda)$ is given by the following explicit formula ($z\in\mathfrak{h}^f$):
\begin{equation*}
     \textnormal{ch}_{H_f(\lambda)}(\tau,z)=q^a \left(\frac{\eta(u\tau)}{\eta(\tau)}\right)^l A(\tau,z)B(\tau,z)C(\tau,z)D(\tau,z),
\end{equation*}
where
\begin{equation*}
    a=\frac{1}{2h^\vee}|u\rho-h^\vee x|^2-\frac{1}{24}\dim\mathfrak{g}_0+\frac{1}{48}\dim\mathfrak{g}_{1/2}+\frac{l(1-u)}{24},
\end{equation*}
\begin{equation*}
    A(\tau,z)=\frac{\prod_{\alpha\in\triangle_+:\alpha(x)>0}(1-q^{\alpha(x)}e^{-2\pi i\alpha(z)})}{\prod_{n=1}^\infty\prod_{\alpha\in\triangle_+:\alpha(x)=1/2}(1-q^{n-\frac{1}{2}}e^{2\pi i\alpha(z)})},
\end{equation*}
\begin{equation*}
    B(\tau,z)=\prod_{n=1}^\infty\prod_{\alpha\in\triangle_+}(1-q^{un+\alpha(x)}e^{-2\pi i\alpha(z)})(1-q^{un-\alpha(x)}e^{2\pi i\alpha(z)}),
\end{equation*}
\begin{equation*}
    C(\tau,z)=\prod_{n=1}^\infty\prod_{\alpha\in\triangle_+:\alpha(x)=0}(1-q^n e^{-2\pi i\alpha(z)})^{-1}(1-q^n e^{2\pi i\alpha(z)})^{-1},
\end{equation*}
\begin{equation*}
    D(\tau,z)=\textnormal{ch}_{\textnormal{L}(\Lambda^0)}(u\tau,-\tau x+z, \frac{\tau}{2u}|x|^2).
\end{equation*}

\underline{Proof}. Since $\mathfrak{g}$ is locally finite on $\textnormal{L}(\lambda)$ (see Remark 6), \cite{Ara15} implies that the normalized Euler-Poincare character of $H_f(\lambda)$ coincides with $\textnormal{ch}_{H_f(\lambda)}$ (see Remark 5). Hence we can apply the character formula (3.3) from
\cite{KRW03},
% where $\hat{R}\textnormal{ch}_\Lambda$ is computed using
 and use formula (5) from \cite{KW17}.

This theorem implies the following corollary.

\underline{Corollary 1}.
(a) If $\lambda\in\hat{Pr}^k$ is as in Theorem 6, then  $H_f(\lambda)$ for a non-zero nilpotent element $f$
% QHR of $\textnormal{L}(\lambda)$, associated to a nilpotent element $f$, 
is zero if and only if $u\leq\theta(x)$. In particular, this holds for $\lambda=k\Lambda_0$.

(b) If $f=f_{min}$, then  $H_f(k\Lambda_0)=W_k^{min}(\mathfrak{g})$, unless $u=1$, when it is 0.

(c) If $f=f_{pr}$ is a principal nilpotent element, then  $H_f(k\Lambda_0)$ is zero iff $u<h$, the Coxeter number of $\mathfrak{g}$.

(d) If $u>\theta(x)$, then $\tilde{W}_k(\mathfrak{g},f)$ is a modular invariant vertex algebra.

\underline{Proof}. (a) Note that $A(\tau,0)$ and $B(\tau,0)$ are non-zero holomorphic functions on the upper half-plane. The same is true for $D(\tau,0)$ since $\textnormal{L}(\Lambda^0)$ is an integrable $\hat{\mathfrak{g}}$-module, hence its character is holomorphic in $z\in\mathfrak{h}$ and $\tau$ in the upper half-plane, and its specialization is non-zero. Finally, $B(\tau,0)$ is zero iff $\alpha(x)=un$ for some $n\in\mathbb{Z}_{\geq1}$. But this happens iff $u\leq\theta(x)$. Indeed, if this inequality holds, then $u=\alpha(x)$ for some $\alpha\in\triangle_+$. Conversely, if $u=\alpha(x)$ for $\alpha\in\triangle_+$, then $u\leq\theta(x)$.
(b) follows from \cite{Ara05}.
(c) holds as well, since $\theta(x)=h-1$ for $f=f_{pr}$.
The proof of (d) is the same as in Example 1.

Corollary 1(a)-(c) is a special case of Theorem 5.16(i) in \cite{Ara15}.

\section{Corrections to previous papers}
We are taking this opportunity to make corrections to our previous papers.
\begin{enumerate}
    \item In the paper \cite{KW17}, Proposition 1a, the condition $gcd(u,h^\vee)=1$ should be removed.
    \item In the paper \cite{KRW03} Theorem 3.2 should be corrected by replacing "the $\mathfrak{g}$-module $M$ is not locally nilpotent with respect to all root spaces $\mathfrak{g}_{-\alpha}$" by "$ch_M$ has a pole at all hyperplanes $T_\alpha$". It is because the "only if" part of Lemma 3.2 is false, hence the reference to it in the proof of Theorem 3.2 should be removed.

    Note that conditions (i)-(iii) on $\alpha$ in Theorem 3.2 mean that $\alpha=-\beta+\beta(x)\delta$, where $\beta\in\triangle_+$, $\beta(x)\in\mathbb{Z}_{>0}$ and $\beta|_{\mathfrak{h}^f}=0$.
\end{enumerate}

%\newpage
%\printbibliography
Acknowledgements. We are grateful to Tomoyuki Arakawa and Chongying Dong
for very important correspondence.


\begin{thebibliography}{99}

\bibitem[Ara05]{Ara05}
T. Arakawa. Representation theory of superconformal algebras
and the Kac-Roan-Wakimoto conjecture,
\emph{Duke Math. J.} 130.3 (2005), pp. 435–478.

\bibitem[Ara15]{Ara15}
T. Arakawa. Associated varieties of modules over
Kac-Moody algebras and $C_2$-cofiniteness of
$W$-algebras, \emph{Int. Math. Res. Not. IMRN} 22
(2015), pp. 11605–11666.

\bibitem[ACK24]{ACK24}
T. Arakawa, T. Creutzig, and K. Kawasetsu. On lisse
non-admissible minimal and principal $W$-algebras, preprint 2024.

\bibitem[AM18]{AM18}
T. Arakawa and A. Moreau. Joseph ideals and lisse minimal
$W$-algebras, \emph{J. Inst. Math. Jussieu} 17.2
(2018), pp. 397–417.

\bibitem[BKK24]{BKK24}
R. Bezrukavnikov, V. Kac, and V. Krylov. Subregular
nilpotent orbits and explicit character formulas for modules over
affine Lie algebras, \emph{Pure Appl. Math. Q.} 20.1
(2024), pp. 81–138.

\bibitem[DLM98]{DLM98}
C. Dong, H. Li, and G. Mason. Twisted representations of
vertex operator algebras, \emph{Math. Ann.} 310.3
(1998), pp. 571–600.

\bibitem[DR18]{DR18}
C. Dong and L. Ren. Congruence property in orbifold theory,
 \emph{Proc. Amer. Math. Soc.} 146.2 (2018), pp. 497–506.

\bibitem[Kac90]{Kac90}
V. G. Kac. Infinite-dimensional Lie algebras. Third,
ed. 1990, pp. xxii+400.

\bibitem[KP84]{KP84}
V. G. Kac and D. H. Peterson. Infinite-dimensional Lie
algebras, theta functions and modular forms, \emph{Adv.
in Math.} 53.2 (1984), pp. 125–264.

\bibitem[KRW03]{KRW03}
V. G. Kac, S.-S. Roan, and M. Wakimoto. Quantum reduction
for affine superalgebras, \emph{Comm. Math. Phys.}
241.2-3 (2003), pp. 307–342.

\bibitem[KW88]{KW88}
V. G. Kac and M. Wakimoto. Modular invariant
representations of infinite-dimensional Lie algebras and
superalgebras, \emph{Proc. Nat. Acad. Sci. U.S.A.} 85.14
(1988), pp. 4956–4960.

\bibitem[KW89]{KW89}
V. G. Kac and M. Wakimoto. Classification of modular
invariant representations of affine algebras, in \emph{ Adv. Ser.
Math. Phys.}, vol 7, 1989, pp. 138–177.

\bibitem[KW04]{KW04}
V. G. Kac and M. Wakimoto. Quantum reduction and
representation theory of superconformal algebras. 
\emph{Adv. Math.} 185.2 (2004), pp. 400–458.

\bibitem[KW14]{KW14}
V. G. Kac and M. Wakimoto. Representations of affine
superalgebras and mock theta functions,
\emph{Transform. Groups} 19.2 (2014), pp. 383–455.

\bibitem[KW17]{KW17}
V. G. Kac and M. Wakimoto. A remark on boundary level
admissible representations. \emph{C. R. Math. Acad. Sci.
Paris} 355.2 (2017), pp. 128–132.

\bibitem[KW18]{KW18}
V. G. Kac and M. Wakimoto. On characters of irreducible
highest weight modules of negative integer level over affine
Lie algebras, in \emph{ Progr. Math.} 2018, vol 326, pp. 235–252.

\bibitem[Kaw18]{Kaw18}
K. Kawasetsu. $W$-algebras with
non-admissible levels and the Deligne exceptional series,
\emph{Int. Math. Res. Not. IMRN} 3 (2018), pp. 641–676.

\bibitem[vEke13]{vEke13}
J. van Ekeren. Modular invariance for twisted modules over a
vertex operator superalgebra, \emph{Comm. Math. Phys.}
322.2 (2013), pp. 333–371.

\bibitem[Wak24]{Wak24}
M. Wakimoto. Vanishing of the quantum reduction of
the Deligne exceptional series representations of negative integer
level,  \href {https://arxiv.org/abs/2408.16560}
{\nolinkurl arXiv:{2408.16560} \texttt{[math.RT]}}.

\bibitem[Zhu96]{Zhu96}
Y. Zhu. Modular invariance of characters of vertex operator
algebras, \emph{J. Amer. Math. Soc.} 9.1 (1996),
pp. 237–302.
\end{thebibliography}
\end{document}